\documentclass[11pt]{article}
\usepackage{mathrsfs}
 \makeatother
\usepackage{amsfonts,amsmath,amssymb}

\title{Existence of Multiple Vortices in\\ Supersymmetric Gauge Field Theory}
\author{Shouxin Chen\\Institute of Contemporary Mathematics\\College of Mathematics and Information Science\\Henan University\\Kaifeng, Henan 475001, PR China\\
\\Yisong Yang\\Department of Mathematics\\Polytechnic Institute of New York University\\Brooklyn, New York 11201, USA}
\date{}

\newtheorem{oldtheorem}{Theorem}[section]
\newtheorem{oldassertion}[oldtheorem]{Assertion}
\newtheorem{oldproposition}[oldtheorem]{Proposition}
\newtheorem{oldremark}[oldtheorem]{Remark}
\newtheorem{oldlemma}[oldtheorem]{Lemma}
\newtheorem{olddefinition}[oldtheorem]{Definition}
\newtheorem{oldclaim}[oldtheorem]{Claim}
\newtheorem{oldcorollary}[oldtheorem]{Corollary}

\newenvironment{theorem}{\begin{oldtheorem}$\!\!\!${\bf.}}{\end{oldtheorem}}

\newcommand{\ddd}{\mbox{d}}
\newcommand{\bea}{\begin{eqnarray}}
\newcommand{\eea}{\end{eqnarray}}
\newcommand{\e}{\mbox{e}}
\newcommand{\pa}{\partial}\newcommand{\Om}{\Omega}

\newcommand{\nn}{\nonumber}

\newcommand{\lbl}[1]{\label{#1}}

\newcommand{\be}{\begin{equation}}
\newcommand{\ee}{\end{equation}}
\newcommand\bes{\begin{eqnarray}}
\newcommand\ees{\end{eqnarray}}
\newcommand{\bess}{\begin{eqnarray*}}
\newcommand{\eess}{\end{eqnarray*}}

\newcommand{\bfR}{\mathbb{R}}

\newcommand\ii{\mbox{i}}

\begin{document}
 \pagestyle{myheadings}

 \thispagestyle{empty}

\setlength{\baselineskip}{16pt}{\setlength\arraycolsep{2pt}

\maketitle

\begin{abstract}
Two sharp existence and uniqueness theorems are presented for solutions of
multiple vortices arising in a six-dimensional brane-world supersymmetric gauge field theory under the general gauge symmetry group $G=U(1)\times SU(N)$
and with $N$ Higgs scalar fields in the fundamental representation of $G$. 
Specifically, when the space of extra dimension is compact so that vortices are hosted in a 2-torus of volume $|\Om|$, the existence of a unique multiple vortex solution 
representing $n_1,\cdots,n_N$ respectively prescribed vortices arising in the $N$ species of the Higgs fields is
established under the explicitly stated necessary and sufficient condition
\[
n_i<\frac{g^2v^2}{8\pi N}|\Om|+\frac1N\left(1-\frac1{N}\left[\frac ge\right]^2\right)n,\quad i=1,\cdots,N,
\]
where $e$ and $g$ are the $U(1)$ electromagnetic and $SU(N)$ chromatic coupling constants,  $v$ measures the energy scale of broken symmetry,
and $n=\sum_{i=1}^N n_i$ is the total vortex number; when the space of extra dimension is the full plane,
the existence and uniqueness of an arbitrarily prescribed $n$-vortex solution of finite energy is always ensured.
These vortices are governed by
a system of nonlinear elliptic equations, which may be reformulated to allow a variational structure.
Proofs of existence are then developed using the methods of calculus of variations.
\end{abstract}

\section{Introduction}
\setcounter{equation}{0}

The concept of solitons is important in quantum field theory. These static solutions, categorized into domain walls, vortices, monopoles, and instantons
and often of topological origins \cite{MS,R,Ryder}, of gauge field equations, give rise to locally concentrated field
configurations and are essential for the description of various fundamental interactions and phenomenologies. Vortices arise in two spatial dimensions
and were first discovered by Abrikosov \cite{Ab} in the form of a mixed state in a type-II superconductor in which the vortex-lines represent partial penetration of magnetic field into
the superconductor as a consequence of partial destruction of superconductivity, in the context of the Ginzburg--Landau theory \cite{GL}.
In quantum field theory, Nielsen and Olesen \cite{NO} showed that vortices arise in the Abelian Higgs model which may be used to model
the so-called dual strings, better known as the Nambu--Goto \cite{Goto} strings, 
which is a basic
construct in string theory \cite{Z}. As a consequence of the Julia--Zee theorem \cite{JZ,SY}, the vortex equations of the Abelian Higgs model 
in the static limit are exactly the Ginzburg--Landau equations. Although these equations are fundamentally important, relatively thorough understanding of their solutions
has only been achieved in a few extreme situations where one assumes either that magnetic field is absent \cite{BBH,E,Lin,Neu}, that the solutions are radially symmetric \cite{BC,CGSY,Pl}, or that a critical coupling is
maintained so that the interaction
between vortices vanishes \cite{JT,WY}. In literature, this last situation is commonly referred to as the self-dual or BPS limit after the pioneering work of Bogomol'nyi \cite{Bo} and 
Prasad--Sommerfield \cite{PS}. In fact, in the area of the study of non-Abelian gauge field equations, it is only in the BPS limit that tractable opportunities are available for
gaining some fair understanding of the solutions of various equations, due to the difficulties associated with the presence of non-Abelian symmetry groups.
The first array of existence results for non-Abelian vortices were obtained  \cite{SYws1,SYws2} for the equations governing electroweak vortices formulated by Ambjorn and Olesen \cite{AO1,AO2,AO3,AO4},
which were later sharpened \cite{BT,CL}. Although the BPS vortices are present only when the coupling constants satisfy a certain critical condition
(in the Ginzburg--Landau theory, this is the interface between type-I and type-II superconductivity; in the Abelian Higgs model, this is when
the masses of the gauge and Higgs bosons coincide), such solutions exhibit a full range of elegant and unambiguous features including exact topological characterization, quantization of flux and energy,
and energy concentration, as anticipated from experimental facts. In recent years, the conceptual power of the BPS vortices in theoretical physics has been particularly witnessed in 
supersymmetric gauge field theory starting with the work of Seiberg and Witten \cite{SW} in an attempt to use non-Abelian color-charged monopoles and vortices to interpret quark confinement
\cite{Man1,Man2, Nambu,tH1,tH2,tH3}. For surveys on this exciting topic, see \cite{Gr,Kon,Shif1,Shif2,Shif3}. Inspired by the importance of non-Abelian vortices in the linear confinement
mechanism through a so-called dual Meissner effect, extending the classical Meissner effect in superconductors, some systematic research has recently been carried out aimed at understanding
the various BPS vortex equations obtained in \cite{EFN,Eto-survey,EI,GJK,ShY2004,ShY-vortex,Shif2,Shif3} and a series of sharp existence theorems have been established \cite{LiebY,LY1,LY2}.
Unlike the BPS vortex equations in the electroweak theory \cite{AO1,AO2,AO3,AO4,BT,SYws1,SYws2}, the elegant structure of the BPS vortex equations in these supersymmetric models allow a complete
understanding of their solutions, despite of the apparent sophistication associated with various underlying non-Abelian gauge groups. The main contribution of the present paper is to prove
two  sharp existence theorems for the BPS vortex equations arising in the supersymmetric $U(1)\times SU(N)$ gauge theory discovered in the work of Eto, Nitta, and Sakai \cite{ENS} in
the context of a six-dimensional brane-world scenario formalism \cite{HLP,HP,Z}.

An outline of the rest of the content of this article is as follows. In the next section, we review the non-Abelian multiple vortex equations of Eto, Nitta, and Sakai \cite{ENS} and state  our
main existence theorems. In Section 3, we describe the nonlinear elliptic equation problem to be studied which is equivalent to the solution problem of the non-Abelian vortex equations.
In Section 4, we consider multiple vortex solutions over a doubly periodic domain and identify a family of necessary conditions for the existence of such solutions. We then prove that these
necessary conditions are also sufficient for existence.  While doing this, we reveal various fine structures of the problem which will be seen to be useful for our later study of planar solutions.
In Section 5, we prove the existence and uniqueness of a weak solution for the vortex equations over the full plane by the variational approach first developed for the scalar (Abelian) situation
\cite{JT}. In Section 6, we first establish pointwise decay properties of solutions near infinity. We then strengthen these results by obtaining some exponential decay estimates. In  Section 7,
we use the obtained exponential decay properties of the solutions to establish various anticipated flux quantization formulas.

\section{Non-Abelian BPS vortex equations and main results}
\setcounter{equation}{0}

In this section, we begin by a review of the non-Abelian vortex equations derived by Eto, Nitta, and Sakai in \cite{ENS}. 
For convenience, we use $x^1$ and $x^2$ to denote the coordinates of extra dimensions in their six-dimensional gauge field theory. We will be sketchy and aimed
at fixing notation since details can be found in  \cite{ENS}. 

We use $\{t^a\}$ to denote the generators of $SU(N)$ ($a=1,\cdots,N^2-1$). Thus any element in the Lie algebra of $SU(N)$, denoted by $su(N)$, may be written as
\be 
\hat{X}=\sum_{a=1}^{N^2-1} X^a t^a.
\ee
We use $\{q_i\}$ ($i=1,\cdots,N$) to denote $N$ hypermultiplets in the fundamental representation of $SU(N)$, i.e., each $q_i$ is a $\mathbb{C}^N$-valued scalar Higgs field
which is such that each of its value is taken to be a column vector.
With $(x^\ell)=(x^1,x^2)$, the gauge-covariant derivatives are given by
\be 
{\cal D}_\ell=\partial_\ell -\ii A_\ell -\ii \hat{A}_\ell,\quad \ell=1,2,
\ee
where $A_\ell$ and $\hat{A}_\ell$ lie in the Lie algebras of $U(1)$ and $SU(N)$, respectively, $\ii=\sqrt{-1}$, so that the induced gauge field strength tensors are defined by
\be 
F_{\ell\ell'}=\partial_\ell A_{\ell'}-\partial_{\ell'} A_\ell,\quad \hat{F}_{\ell\ell'}=\partial_\ell \hat{A}_{\ell'}-\partial_{\ell'}\hat{A}_\ell-\ii[\hat{A}_\ell,\hat{A}_{\ell'}],\quad \ell,\ell'=1,2,
\ee
where $[\cdot,\cdot]$ denotes the matrix commutator.

With the above notation, the non-Abelian vortex equations derived in \cite{ENS} (see also \cite{Au,HT}) are
\bea
F^a_{12}&=& \frac{g^2}2 \sum_{i=1}^{N}q_i^\dagger t^a q_i,\quad a=1,\cdots,N^2-1,\label{2.4}\\
F_{12}&=&\frac{e^2}2 \left(\sum_{i=1}^{N}q^\dagger_i q_i-v^2\right),\label{2.5}\\
{\cal D}_1 q_i&=&-{\cal D}_2 q_i,\quad i=1,\cdots,N,\label{2.6}
\eea
where $e,g,v>0$ are coupling parameters for which $e$ represents Abelian (electromagnetic) gauge coupling, $g$ non-Abelian (nuclear) gauge coupling,  $v$ the energy
scale of the spontaneously broken ground state (vacuum), and $^\dagger$ denotes the Hermitian conjugate. 

Now, following \cite{ENS}, we collectively rewrite the Higgs fields and gauge fields in forms of $N\times N$ matrices,
\be 
q=(q_1,\cdots,q_N),\quad \check{A}_\ell=A_\ell I_N+\hat{A}_\ell,\quad\ell=1,2,
\ee
where $I_N$ denotes the $N\times N$ identity matrix. Then the ansatz
\bea 
q&=&\frac1{\sqrt{N}}\mbox{diag}\{\phi_1,\cdots,\phi_N\},\\
\check{A}_\ell&=&\mbox{diag}\{B^1_\ell,\cdots,B^N_\ell\},\quad \ell=1,2,
\eea
where $\phi_i$ are complex-valued scalar fields and $B^i=(B^i_\ell)$ are real-valued vector fields, $i=1,\cdots,N$, further reduces the equations
(\ref{2.4})--(\ref{2.6}) into \cite{ENS}:
\bea
4B^i_{12}&=&\frac{g^2}N(v^2-|\phi_i|^2)+\left(e^2-\frac{g^2}N\right)\left(v^2-\frac1N\sum_{i=1}^N|\phi_i|^2\right),\label{2.10}\\
(\pa_1+\ii\pa_2)\phi_i&=&\ii(B^i_1+\ii B^i_2)\phi_i,\label{2.11}
\eea
for $i=1,\cdots,N$.

In view of \cite{JT}, we see from (\ref{2.11}) that the zeros of each $\phi_i$ are isolated with integer multiplicities. We may use
$Z(\phi_i)$ to denote the set of zeros of $\phi_i$,
\be \label{2.12}
Z(\phi_i)=\{p_{i,1},\cdots,p_{i,n_i}\},\quad i=1,\cdots, N,
\ee
so that the repetitions among the points $p_{i,s}$ ($s=1,\cdots,n_i$) take account of the multiplicities of these zeros. The first term on the right-hand side of (\ref{2.10}) indicates that
these zeros enhance the `vorticity' field $B^i_{12}$ ($i=1,\cdots,N$). However, the other terms on the right-hand side of (\ref{2.10}) complicate the situation so that we are not able to
assert that the maxima of $B^i_{12}$ are achieved at $Z(\phi_i)$ ($i=1,\cdots,N$), which is what makes the problem interesting and challenging.

There are two partial differential equation problems to be studied. 

The first one concerns the solutions of (\ref{2.10})--(\ref{2.11}) with prescribed zero sets given in (\ref{2.12}) over a bounded
domain $\Om$ in $\bfR^2$ so that the field configurations are periodic modulo the 't Hooft boundary condition \cite{tH,WY} at the boundary of $\Om$. 
The physical relevance of this problem is that such solutions may arise when the extra dimensions are compactified as two-dimensional tori. For this problem, here is our result.

\begin{theorem}\label{theorem2.1}
Consider the non-Abelian BPS vortex equations (\ref{2.10})--(\ref{2.11}) for the field configurations $(\phi_i,B^i_\ell)$ ($i=1,\cdots,N$) 
over a doubly periodic domain $\Om$ with the given sets of zeros stated in (\ref{2.12}) so that $\phi_i$ has exactly $n_i$ zeros in $\Om$ ($i=1,\cdots,N$).
Then a solution exists if and only if the  conditions
\be \label{NS}
n_i<\frac{g^2v^2}{8\pi N}|\Om|+\frac1N\left(1-\frac1{N}\left[\frac ge\right]^2\right)n,\quad i=1,\cdots,N,
\ee 
are fulfilled simultaneously, where $n=\sum_{i=1}^N n_i$ is the total number of the $N$ species of vortices. Furthermore, if a solution exists, it must be uniquely determined by the sets of zeros given in (\ref{2.12}), up to gauge transformations.
\end{theorem}

Although the condition (\ref{NS}) seems complicated, we may sum up the vortex numbers, $n_1,\cdots,n_N$, to get a simple consequence of the inequalities stated in (\ref{NS}), in the form
\be \label{NS1}
n<\frac{e^2 v^2}{8\pi}N|\Om|.
\ee
It is interesting to notice that, now, the non-Abelian coupling constant $g$ does not enter the condition (\ref{NS1}) but the integer $N$ which measures the `size' of the non-Abelian symmetry.

The second problem concerns the solutions of (\ref{2.10})--(\ref{2.11}) over the full plane $\bfR^2$. The form of the energy density described in \cite{ENS} and requirement of finite energy
naturally impose the boundary condition
\be \label{BC}
\lim_{|x|\to\infty}|\phi_i|(x)=v,\quad i=1,\cdots,N,
\ee 
for solutions. For this problem, here is our result.

\begin{theorem}\label{theorem2.2}
The non-Abelian BPS vortex equations (\ref{2.10})--(\ref{2.11}) over the full plane $\bfR^2$ described by the field configurations $(\phi_i,B^i_\ell)$ ($i=1,\cdots,N$) 
subject to the given sets of zeros stated in (\ref{2.12}) and the boundary condition (\ref{BC}) always have a unique solution. Moreover, such a solution realizes the boundary
condition (\ref{BC}) exponentially fast. More precisely, for an arbitrarily small number $\varepsilon\in (0,1)$, there hold
\be 
\left||\phi_i|-v\right| +\sum_{\ell=1}^2\left|(\partial_\ell -\mbox{\rm i} B^i_\ell)\phi_i\right|^2 +\left| B^i_{12}\right|=\mbox{\rm O}(\mbox{\rm e}^{-(1-\varepsilon)\frac g{\sqrt{N}}|x|}),\quad i=1,\cdots,N,
\ee
for $|x|$ sufficiently large.
\end{theorem}

We note that, intuitively, Theorem \ref{theorem2.2} may be reinterpreted in view of Theorem \ref{theorem2.1} in the context of vortices in a domain of infinite volume. In such a situation,
the condition (\ref{NS}) of course becomes superfluous.

\begin{theorem}\label{theorem2.3}
In both situations stated as in Theorems \ref{theorem2.1} and \ref{theorem2.2}, the fluxes are quantized quantities given by the expressions
\be \label{Bi}
\int B^i_{12}\,{\mbox{\rm d}} x=2\pi n_i,\quad i=1,\cdots,N,
\ee
which are seen to be determined by the algebraic numbers of zeros of the complex scalar fields $\phi_i$ ($i=1,\cdots,N$), with the integral evaluated over
the doubly periodic domain $\Omega$ or $\mathbb{R}^2$, respectively.
\end{theorem}

This result suggests that, although is not clear whether the vorticity field $B^i_{12}$ peaks at the zeros of the order parameter $\phi_i$, the total number of zeros of
$\phi_i$ determines the total vorticity (or flux) generated from the field $B^i_\ell$ in the full domain. For this reason, we may still regard the locations of the zeros of $\phi_i$ 
as the centers of vortices and refer to $n_i$ as the $i$th vortex number, which extends the concept of vortices in the classical Abelian Higgs model \cite{JT}. As a consequence,
the integer $n$ defined in (\ref{NS}) is well justified to be called the total vortex number of the solution.

With the quantized fluxes given in (\ref{Bi}), the total non-Abelian vortex (string) tension, $T_{\mbox{NA}}$, may be computed  \cite{ENS} to assume the elegant exact value
\be 
T_{\mbox{NA}}=\frac{v^2}N\sum_{i=1}^N\int B^i_{12}\,\ddd x=\frac{2\pi v^2}N\sum_{i=1}^N n_i.
\ee
 
The above theorems will be established in the subsequent sections.

\section{System of nonlinear elliptic equations}
\setcounter{equation}{0}

To proceed, we now adapt the complexified variables and derivatives defined by
\be 
z=x^1+\ii x^2,\quad B^i=B_1^i+\ii B^i_2,\quad \pa=\frac12(\pa_1-\ii\pa_2),\quad \overline{\pa}=\frac12(\pa_1+\ii\pa_2),
\ee
and convert the equations (\ref{2.10})--(\ref{2.11}) into the following system
\bea 
4\ii(\pa B^i-\overline{\pa}\overline{B}^i)&=&\frac{g^2}N(|\phi_i|^2-v^2)+\left(e^2-\frac{g^2}N\right)\left(\frac1N\sum_{i=1}^N|\phi_i|^2-v^2\right),\label{2.13}\\
\overline{\pa}\ln\phi_i&=&\frac\ii2 B^i,\label{2.14}
\eea
away from the possible zeros of $\phi_i$, for $i=1,\cdots,N$.

Inserting (\ref{2.14}) or $B^i=-2\ii\overline{\pa}\ln\phi_i$  into (\ref{2.13}) ($i=1,\cdots,N$) and using the
relation $\Delta=\pa_1^2+\pa_2^2=4\pa\overline{\pa}=4\overline{\pa}\pa$, we arrive at the equations
\be
\Delta \ln|\phi_i|^2=\frac{g^2}{2N}(|\phi_i|^2-v^2)+\frac12\left(e^2-\frac{g^2}N\right)\left(\frac1N\sum_{i=1}^N|\phi_i|^2-v^2\right),
\ee
away from the zeros of $\phi_i$, $i=1,\cdots,N$. Thus, with the notation in (\ref{2.12}) and the new variables
\be 
u_i=\ln|\phi_i|^2,\quad i=1,\cdots,N,
\ee
we obtain the following system of nonlinear elliptic equations 
\be\label{S}
\Delta u_i=\frac{g^2}{2N}(\e^{u_i}-v^2)+\frac12\left(e^2-\frac{g^2}N\right)\left(\frac1N\sum_{i=1}^N\e^{u_i}-v^2\right)+4\pi\sum_{s=1}^{n_i}\delta_{p_{i,s}}(x),\quad i=1,\cdots,N,
\ee
governing the interaction of $N$ species of vortices located at the prescribed set of points
\be 
Z=\cup_{i=1}^N Z(\phi_i),
\ee
which are the set of zeros of the complex scalar fields $\phi_1,\cdots,\phi_N$.

We shall consider the solutions of the system (\ref{S}) over a doubly periodic domain (a 2-torus) $\Om$ realized by the 't Hooft periodic boundary condition
\cite{tH,WY} and over the full plane $\bfR^2$. In the latter situation, we need to observe the boundary condition (\ref{BC}) at infinity. That is,
\be \label{bc}
\lim_{|x|\to\infty} u_i=2\ln v,\quad i=1,\cdots,N.
\ee

In the next section, we shall first concentrate on the doubly periodic situation.

\section{Necessary and sufficient condition for doubly periodic solutions}
\setcounter{equation}{0}

In this section, we study the  equations (\ref{S}) defined over a
doubly periodic domain $\Omega$. We conveniently rewrite these equations as
  \bes
\Delta u_i=\sum_{j=1}^{N}a_{ij}(\e^{u_j}-v^2)+4\pi\sum_{s=1}^{n_i}\delta_{p_{i,s}}(x), \ \ \
i=1,\cdots,N,
 \lbl{3.1}\ees
where
   \bes
a_{ij}=\frac{1}{N}\left(\frac{e^2}{2}-\frac{g^2}{2N}\right)+\delta_{ij}\frac{g^2}{2N},\ \
i,j=1,\cdots,N.
   \lbl{3.2}\ees

Let $u^0_i$ be a solution to
   \bes
\Delta u_i^0=-\frac{4\pi n_i}{|\Omega|}+4\pi\sum_{s=1}^{n_i}\delta_{p_{i,s}}(x), \ \ \
i=1,\cdots,N.
 \lbl{3.3}\ees
(cf.  \cite{Aubin}.)
The substitutions
   \bess
u_i=u_i^0+U_i,  \ \ \ i=1,\cdots,N,
    \eess
recast the equation (\ref{3.1}) into
  \bes
\Delta U_i=\sum_{j=1}^{N}a_{ij}(\e^{u_j^0+U_j}-v^2)+\frac{4\pi n_i}{|\Omega|}, \ \ \
i=1,\cdots,N.
 \lbl{3.4}\ees

We use boldfaced letters to denote column vectors in $\bfR^n$. Thus, we set
   \bess
&& \textbf{U}=(U_1,\cdots,U_N)^\tau,\ \ \ \textbf{G}=\left(\e^{u_1^0+U_1},\cdots,\e^{u_N^0+U_N}\right)^\tau,\\
&& \textbf{F}=\left(\frac{4\pi n_1}{|\Omega|}-v^2\sum_{j=1}^{N}a_{1j},
\cdots,\frac{4\pi n_N}{|\Omega|}-v^2\sum_{j=1}^{N}a_{Nj}\right)^\tau
\equiv(f_1,\cdots,f_N)^\tau,
   \eess
and let $A=(a_{ij})_{N\times N}$ be the $N\times N$ matrix defined by (\ref{3.2}). Then
  \bess
A=\left(\begin{array}{ccccc}
 a+b &\ a &\ a &\cdots &\  a \\[1mm]
 a &\ a+b &\ a &\cdots &\  a\\[1mm]
 \cdots &\ \cdots &\ \cdots  &\ \cdots &\  \cdots\\[1mm]
 a &\ a &\ a &\cdots &\  a+b
\end{array}\right), \nonumber
   \eess
where $a=\frac{1}{N}\left(\frac{e^2}{2}-\frac{g^2}{2N}\right)$ and $b=\frac{g^2}{2N}$ which should not be confused with the group indices used in Section 2.

Now the equations (\ref{3.4}) can be written in the vector form
   \bes
\Delta \textbf{U}=A \textbf{G}+\textbf{F}.
   \lbl{3.5}\ees
This system looks difficult to approach. In order to tackle it, we shall seek for a variational principle.

To find a variational principle, we need to use the property of the
matrix $A$. It is easy to check that the matrix $A$ is positive
definite and its eigenvalues are
\be  \label{lambda}
\lambda_1=Na+b=\frac{e^2}{2},\quad \lambda_2=\cdots=\lambda_N=b=\frac{g^2}{2N}. 
\ee
Then, by the Cholesky decomposition theorem \cite{GO}, we know that
there is a unique upper triangular $N\times N$ matrix $T=(t_{ij})$
for which all the diagonal entries are positive, i.e., $t_{ii}>0$,
$i=1,\cdots,N$, such that
   \bes
A=T^\tau T.
   \lbl{3.6}\ees
   
In fact, by direct computation, we have
   \bess
&&t_{11}=\sqrt{a+b},\quad
t_{12}=t_{13}=\cdots=t_{1N}=\frac{a}{t_{11}}\equiv\alpha_1>0,\\
&&t_{22}=\sqrt{(a+b)-\alpha^2_1},\quad
 t_{23}=t_{24}=\cdots=t_{2N}=\frac{a-\alpha_1^2}{t_{22}}\equiv\alpha_2>0,\\
&&\cdots\cdots\cdots\\
&&t_{(N-1)(N-1)}=\sqrt{(a+b)-\sum_{i=1}^{N-2}\alpha_i^2},\quad
 t_{(N-1)N}=\frac{a-\sum_{i=1}^{N-2}\alpha_i^2}{t_{(N-1)(N-1)}}\equiv\alpha_{N-1}>0,\\
&&t_{NN}=\sqrt{(a+b)-\sum_{i=1}^{N-1}\alpha_i^2}.
   \eess

Set $\textbf{v}=(v_1,\cdots,v_N)^\tau$, $L=(T^\tau)^{-1}\equiv(l_{ij})_{N\times N}$.
We introduce the new variable vector
   \bes
\textbf{v}=L\textbf{U} \ \  {\rm or}\ \ \textbf{U}=L^{-1}\textbf{v}=T^\tau{\bf v}.
   \lbl{3.7}\ees
Then (\ref{3.5}) takes the form
    \bes
\Delta \textbf{v}=T\textbf{G}+L\textbf{F}.
   \lbl{3.8}\ees
   
With the convention $\alpha_N=0$, we may write (\ref{3.8}) in the component form
    \bes
\Delta v_i=t_{ii}\e^{u_i^0+t_{ii}v_i+\sum_{k=1}^{i-1}\alpha_kv_k}+\alpha_i\sum_{j=i+1}^N
\e^{u_j^0+t_{jj}v_j+\sum_{k=1}^{j-1}\alpha_kv_k}+\sum_{j=1}^il_{ij}f_j, \ \ \ i=1,\cdots,N.
   \lbl{3.9}\ees
   
It is easy to check that the above system of equations (\ref{3.9}) are
the Euler--Lagrange equations of the functional
   \bes
I(\textbf{v})=\int_\Omega\left\{\frac{1}{2}\sum_{i=1}^N|\nabla v_i|^2+
\sum_{i=1}^N\e^{u_i^0+t_{ii}v_i+\sum_{k=1}^{i-1}\alpha_k v_k}
+\sum_{i=1}^N\left(\sum_{j=1}^i l_{ij} f_j\right)v_i\right\}\ddd x.
   \lbl{3.10}\ees

Setting
\be \label{qi}
q_i=\int_\Omega \e^{u_i^0+t_{ii}v_i+\sum_{k=1}^{i-1}\alpha_kv_k}\,\ddd x,\quad i=1,\cdots,N,
\ee 
(note that these $q_i$'s should not be confused with the $N$ scalar ``quark" fields in the fundamental representation of the gauge group
denoted by the same notation in Section 2) we may integrate (\ref{3.9}) to obtain
\be \label{tba}
t_{ii}q_i+\alpha_i\sum_{j=i+1}^Nq_j=-|\Om|\sum_{j=1}^i l_{ij}f_j\equiv p_i,\quad i=1,\cdots,N.
\ee 
Since $t_{ii}>0$ ($i=1,\cdots,N$) and $\alpha_i>0$ ($i=1,\cdots, N-1$), the definition of $q_i$ ($i=1,\cdots,N$) given in (\ref{qi}) and the relation (\ref{tba}) lead to
the necessary condition
\be\label{pi}
p_i>0,\quad i=1,\cdots,N,
\ee
which appears to be complicated.

In order to arrive at an explicit form of the necessary condition,
 recall the structure of the matrix $T$ given in (\ref{3.6}). Thus, with ${\bf p}=(p_1,\cdots,p_N)^\tau$ and ${\bf q}=(q_1,\cdots,q_N)^\tau$, we can rewrite (\ref{tba}) in the matrix form
\be \label{Tq}
T{\bf q}=-|\Om|L{\bf F}={\bf p}.
\ee
With $L=(T^\tau)^{-1}$, we can solve (\ref{Tq}) to get
\be \label{qeq}
{\bf q}=-|\Om|T^{-1}L{\bf F}=-|\Om| A^{-1}{\bf F}.
\ee

On the other hand, since for any invertible $N\times N$ matrix $D$ and the column vectors ${\bf X}$ and ${\bf Y}$ in $\bfR^N$ satisfying
${\bf Y}^\tau D^{-1} {\bf X}\neq1$, the matrix $M=D-{\bf X}{\bf Y}^\tau$ is invertible and
\be \label{M}
M^{-1}=(D-{\bf X}{\bf Y}^\tau)^{-1}=(I+[1-{\bf Y}^\tau D^{-1} {\bf X}]^{-1} D^{-1}{\bf X}{\bf Y}^\tau)D^{-1}.
\ee 
Applying the formula (\ref{M}) to the matrix
\be 
A=\mbox{diag}\{b,\cdots,b\}-(-a,\cdots,-a)^\tau (1,\cdots,1),
\ee
we have 
\be \label{A}
A^{-1}=\frac1{b(Na+b)}\left(\begin{array}{cccc}(N-1)a+b&-a&\cdots&-a\\-a&(N-1)a+b&\cdots&-a\\\cdots&\cdots&\cdots&\cdots\\-a&-a&\cdots&(N-1)a+b\end{array}\right).
\ee 
Inserting (\ref{A}) into (\ref{qeq}),  we obtain the solution
\bea \label{q1}
q_i&=&-|\Om|\left(\frac1bf_i-\frac{a}{b(Na+b)}\sum_{j=1}^N f_j\right)\nn\\
&=&v^2|\Om|+\frac{4\pi a}{b(Na+b)}n-\frac{4\pi}b n_i,\quad i=1,\cdots,N,
\eea
where
\be 
n=\sum_{j=1}^N n_i
\ee
is the total vortex number. Substituting the values of the constants $a$ and $b$ in terms of the coupling constants $e$ and $g$ in (\ref{q1}), we have
\be \label{qi2}
q_i=v^2|\Om|+8\pi\left(\frac1{g^2}-\frac1{N e^2}\right)n-\frac{8\pi N}{g^2} n_i>0,\quad i=1,\cdots,N,
\ee 
which lead us to the necessity of the condition (\ref{NS}).

Below, we shall show that, under the condition (\ref{NS}), the
system (\ref{3.9}) has a solution. We will use a direct minimization
method as in \cite{LiebY}. Furthermore, in order to gain more insight to the technical structure of the problem,
we shall also sketch a constrained minimization method to approach the problem.


We use $W^{1,2}(\Omega)$ to denote the usual Sobolev space of
scalar-valued or vector-valued $\Omega$-periodic $L^2$-functions whose derivatives
are also in $L^2(\Omega)$. In the scalar case, we may decompose $W^{1,2}(\Omega)$
into $W^{1,2}(\Omega)=\mathbb{R}\oplus \dot{W}^{1,2}(\Omega)$ so that any $f\in W^{1,2}(\Omega)$
can be expressed as
   \bes
f=\underline{f}+\dot{f},\quad \underline{f}\in\mathbb{R},\quad \dot{f}\in\dot{W}^{1,2}(\Omega),\quad
\int_\Omega \dot{f}\,\ddd x=0.
   \lbl{3.13}\ees
   
It is useful to recall the Moser--Trudinger inequality \cite{Aubin,Fo}
    \bes
\int_\Omega \e^u \,\ddd x\leq C\exp\left(\frac{1}{16\pi}\int_\Omega|\nabla u|^2\, \ddd x\right),\ \ \ \
u\in\dot{W}^{1,2}(\Omega).
   \lbl{3.14}\ees
With (\ref{3.14}), it is clear that the functional $I$ defined by
(\ref{3.10}) is a $C^1$-functional with respect to its argument
$(v_1,\cdots,v_N)\in W^{1,2}(\Omega)$, which is strictly convex and lower semi-continuous
in terms of the weak topology of $W^{1,2}(\Omega)$.

We can suppress the functional $I$ given in (\ref{3.10})  into the form
   \be\label{3.16}
I(\textbf{v})=\int_\Omega\left(\frac{1}{2}\sum_{i=1}^N|\nabla v_i|^2+
\sum_{i=1}^N\e^{u_i^0+t_{ii}v_i+\sum_{k=1}^{i-1}\alpha_kv_k}\right)\,\ddd x-\sum_{i=1}^N p_i \underline{v}_i.
\ee

Applying the Jensen inequality and using the fact that $t_{ii}>0$ ($i=1,\cdots,N$) and $\alpha_i>0$ ($i=1,\cdots,N-1$), we have
   \bea
\int_\Omega \e^{u_i^0+t_{ii}(\underline{v}_i+\dot{v}_i)+
\sum_{k=1}^{i-1}\alpha_k(\underline{v}_k+\dot{v}_k)}\,\ddd x
&\geq&
|\Omega|\exp\left(-\frac{1}{|\Omega|}
\int_\Omega u_i^0\,\ddd x\right)\exp\left(t_{ii}\underline{v}_i+
\sum_{k=1}^{i-1}\alpha_k\underline{v}_k\right)\nonumber\\
&\equiv&\sigma_i \e^{t_{ii}\underline{v}_i+
\sum_{k=1}^{i-1}\alpha_k\underline{v}_k},\ \ \ i=1,\cdots,N.
   \lbl{3.17}\eea

Substituting (\ref{3.17}) into (\ref{3.16}), we have
   \be
I(\textbf{v})-\frac{1}{2}\int_\Omega\sum_{i=1}^N|\nabla \dot{v}_i|^2\,\ddd x
\geq\sum_{i=1}^N\sigma_i \e^{t_{ii}\underline{v}_i+
\sum_{k=1}^{i-1}\alpha_k\underline{v}_k}-\sum_{i=1}^N p_i\underline{v}_i.
   \lbl{3.18}\ee

In order to control the lower bound of (\ref{3.18}), 
 we recall the relation $T{\bf q}={\bf p}$ with the quantities $q_i$ ($i=1,\cdots,N$) as defined in (\ref{qi2}). That is,
\be 
p_i=t_{ii}q_i+\alpha_i\sum_{j=i+1}^N q_j,\quad i=1,\cdots,N.
\ee
Therefore we have 
\be 
\sum_{i=1}^Np_i \underline{v}_i=\sum_{i=1}^N q_i\left(t_{ii}\underline{v}_i+\sum_{k=1}^{i-1} \alpha_k \underline{v}_k\right).
\ee 
Now set
\be\label{wi}
\underline{w}_i=t_{ii}\underline{v}_i+\sum_{k=1}^{i-1} \alpha_k \underline{v}_k,\quad i=1,\cdots,N.
\ee
Then we arrive at
\be
I(\textbf{v})-\frac{1}{2}\int_\Omega\sum_{i=1}^N|\nabla \dot{v}_i|^2\,\ddd x\geq\sum_{i=1}^N\left(\sigma_i\e^{\underline{w}_i}-q_i \underline{w}_i\right).
   \lbl{x3.21}\ee
Thus, using the elementary inequality
    \bes
\frac{a}{b}\left(1-\ln\left[\frac{a}{bc}\right]\right)\leq ce^{bx}-ax, \quad a,b,c>0,\quad x\in\mathbb{R},
   \lbl{3.20}\ees
    in (\ref{x3.21}), we obtain the coercive lower bound
\be\label{3.21a}
I(\textbf{v})-\frac{1}{2}\int_\Omega\sum_{i=1}^N|\nabla \dot{v}_i|^2\,\ddd x\geq\sum_{i=1}^N q_i\left(1+\ln\left[\frac{\sigma_i}{q_i}\right]\right).
\ee

It follows from (\ref{3.21a}) that $I(\textbf{v})$ is bounded from below and
we may consider the following direct minimization problem
   \bes
\eta\equiv\inf\left\{I(\textbf{v})| \ \textbf{v}\in W^{1,2}(\Omega)\right\}.
   \lbl{3.22}\ees
   
Let $\{(v_1^{(n)},\cdots,v_N^{(n)})\}$ be a minimizing sequence of (\ref{3.22}).
Since the function
   \bes
F(u)=\sigma e^u-q u,
   \lbl{3.23}\ees
where $\sigma,q>0$ are constants, enjoys the property that $F(u)\rightarrow\infty$
as $u\rightarrow\pm\infty$, we see from (\ref{x3.21}) that the sequences $\{\underline{w}_i^{(n)}\}$
$(i=1,\cdots,N)$ are all bounded where $\underline{w}^{(n)}_i$ is defined by (\ref{wi}) by setting
$\underline{w}_i=\underline{w}_i^{(n)}$ and $\underline{v}_i=\underline{v}_i^{(n)}$ ($i=1,\cdots,N, n=1,2,\cdots$). Inverting the transformation (\ref{wi}),
we see that the sequences $\{\underline{v}_i^{(n)}\}$ ($i=1,\cdots,N$) are also bounded.
Without loss of generality, we may assume
    \bes
\underline{v}_i^{(n)}\rightarrow \ {\rm some}\ {\rm point}\ \underline{v}_i^{(\infty)}\in\mathbb{R} \ \
{\rm as}\ n\rightarrow\infty,\ i=1,\cdots,N.
   \lbl{3.24}\ees

In addition, using (\ref{3.21a}), we conclude that $\{\nabla \dot{v}_i^{(n)}\}$
$(i=1,\cdots,N)$ are all bounded in $L^2(\Omega)$. Therefore, it follows from the Poincar$\acute{e}$
inequality that the sequences $\{\dot{v}_i^{(n)}\}$
$(i=1,\cdots,N)$ are bounded in $W^{1,2}(\Omega)$. Without loss of generality, we may assume
    \bes
\dot{v}_i^{(n)}\rightarrow \ {\rm some}\ {\rm element}\ \dot{v}_i^{(\infty)}\in W^{1,2}(\Omega) \ {\rm weakly} \
{\rm as}\ n\rightarrow\infty,\ i=1,\cdots,N.
   \lbl{3.25}\ees
   
Obviously, $\dot{v}_i^{(\infty)}\in \dot{W}^{1,2}(\Omega)$ $(i=1,\cdots,N)$.
Set $v_i^{(\infty)}=\underline{v}_i^{(\infty)}+\dot{v}_i^{(\infty)}$ $(i=1,\cdots,N)$.
Then (\ref{3.24}) and (\ref{3.25}) lead us to see that $v_i^{(n)}\rightarrow v_i^{(\infty)}$ weakly in
$W^{1,2}(\Omega)$ as $n\rightarrow\infty$ $(i=1,\cdots,N)$. The weakly lower semi-continuity
of $I$ enables to conclude that $(v_1^{(\infty)},\cdots,v_N^{(\infty)})$
solves (\ref{3.22}), which is a critical point of $I$ and a classical solution to the system (\ref{3.9}) in view of the standard elliptic theory.

Since the matrix $A$ is positive definite, it is easy to check that the functional $I$ is strictly convex in $W^{1,2}(\Omega)$.
So it has at most one critical point in $W^{1,2}(\Om)$, which establishes the uniqueness of the
solution to the equations (\ref{3.9}).

Below, we will develop our methods further by presenting a constrained minimization approach to the problem.
For convenience, we rewrite the constraints (\ref{qi}) collectively as
   \bes
J_i(\textbf{v})=\int_\Omega \e^{u_i^0+t_{ii}v_i+\sum_{k=1}^{i-1}\alpha_kv_k}\,\ddd x=q_i, \ \
i=1,\cdots,N.
   \lbl{3.26}\ees
Recall that the values of $q_1,\cdots,q_N$ are given by (\ref{qi2}) which are obtained by solving the system of equations (\ref{tba}), or
\be \label{tpq}
t_{ii}q_i+\alpha_i\sum_{j=i+1}^N q_j+|\Om|\sum_{j=1}^i l_{ij}f_j=0,\quad i=1,\cdots,N.
\ee 

We consider the multi-constrained minimization problem
   \bes
\eta\equiv\inf\{I(\textbf{v})| \ \textbf{v}\in W^{1,2}(\Omega), J_1({\bf v})=q_1,\cdots,J_N({\bf v})=q_N\}.
   \lbl{3.27}\ees

Suppose that (\ref{3.27}) allows a solution, say ${\bf v}=({v}_1,\cdots,{v}_N)$.
Then there are numbers (the Lagrange multipliers) in $\mathbb{R}$, say $\xi_1,\cdots,\xi_N$, such that, for $i=1,\cdots,N$,
   \bes
&&\int_\Omega \left\{\nabla v_i\cdot\nabla w_i+\left(t_{ii}\e^{u_i^0+t_{ii}v_i+\sum_{k=1}^{i-1}\alpha_kv_k}+\alpha_i\sum_{j=i+1}^N
\e^{u_j^0+t_{jj}v_j+\sum_{k=1}^{j-1}\alpha_kv_k}+\sum_{j=1}^il_{ij}f_j\right)w_i\right\}\,\ddd x
\nonumber\\
&=&\xi_it_{ii}\int_\Omega \e^{u_i^0+t_{ii}v_i+\sum_{k=1}^{i-1}\alpha_kv_k}w_i\,\ddd x+\alpha_i\sum_{j=i+1}^N
\xi_j\int_\Omega
\e^{u_j^0+t_{jj}v_j+\sum_{k=1}^{j-1}\alpha_kv_k} w_i\,\ddd x,
    \lbl{3.28}\ees
where $w_1,\cdots,w_N$ are test functions. Letting $w_1=\cdots=w_N=1$ in the above equations and applying (\ref{tpq}), we arrive at
   \bes
\xi_it_{ii}q_i+\alpha_i\sum_{j=i+1}^N
\xi_j q_j=0, \ \ i=1,\cdots,N.
  \lbl{3.29} \ees
Consequently,  $\xi_1=\cdots=\xi_N=0$ so that (\ref{3.28}) is exactly the weak form of the system (\ref{3.9}). In other words,
the Lagrange multipliers disappear automatically and a solution of (\ref{3.27})
solves (\ref{3.9}). Hence, it suffices to find a solution to (\ref{3.27}).

In order to approach (\ref{3.27}), we use the notation (\ref{3.13}) to rewrite the constraints (\ref{3.26}) as
\be 
\e^{t_{ii}\underline{v}_i+\sum_{k=1}^{i-1}\alpha_k\underline{v}_k}\int_\Om \e^{u_i^0+t_{ii}\dot{v}_i+\sum_{k=1}^{i-1}\alpha_k\dot{v}_k}\,\ddd x=q_i,\quad
i=1,\cdots,N,
\ee 
which may be resolved to yield
   \bes
\underline{v}_i=\sum_{j=1}^il_{ij}(\ln q_j-\ln J_j(\dot{\bf v})), \ \ i=1,\cdots,N,
  \lbl{3.31}\ees
where $\dot{\bf v}=(\dot{v}_1,\cdots,\dot{v}_N)^\tau$ and 
$l_{ij}$ are the entries of the lower triangular matrix $L=(T^\tau)^{-1}$ with $T=(t_{ij})$.

To proceed, we use the constraints (\ref{3.26}) to rewrite the action functional (\ref{3.10}) or (\ref{3.16}) as
\bea \label{ex1}
I({\bf v})-\frac12\sum_{i=1}^N\int_\Om|\nabla \dot{v}_i|^2\,\ddd x&=&\sum_{i=1}^N q_i-\sum_{i=1}^N p_i\underline{v}_i\nn\\
&=&\sum_{i=1}^N\sum_{j=1}^i p_i l_{ij} \ln J_j(\dot{\bf v})+\sum_{i=1}^N\left(q_i-p_i\sum_{j=1}^i  l_{ij}\ln q_j\right),
\eea
where we have inserted (\ref{3.31}). However, using the relation ${\bf q}=L^\tau {\bf p}$, we can rewrite (\ref{ex1}) as
  \be\label{ex2}
I({\bf v})-\frac12\sum_{i=1}^N\int_\Om|\nabla \dot{v}_i|^2\,\ddd x=\sum_{i=1}^N q_i \ln J_i(\dot{\bf v})-C_0,
\ee
where $C_0$ is a constant depending only on $L$, ${\bf p}$, and $\bf q$. By virtue of the Jensen inequality, we have
\be \label{ex3}
J_i(\dot{\bf v})\geq|\Om|\exp\left(\int_\Om u_i^0\,\ddd x\right)\equiv \sigma_i,\quad i=1,\cdots,N.
\ee 
Using the condition $q_1,\cdots,q_N>0$ and the lower bound (\ref{ex3}) in (\ref{ex2}), we obtain the following coercive inequality
 \be\label{ex4}
I({\bf v})-\frac12\sum_{i=1}^N\int_\Om|\nabla \dot{v}_i|^2\,\ddd x\geq\sum_{i=1}^N q_i \ln \sigma_i-C_0.
\ee

 Now the proof of solvability of (\ref{3.27}) follows from a standard argument. 

In fact, let
$\{(v_1^{(n)},\cdots,v_N^{(n)})\}$ be a minimizing sequence of (\ref{3.27}).
 In view of (\ref{ex4}) and the
Poincar$\acute{e}$ inequality, we see that $\{(\dot{v}_1^{(n)},\cdots,\dot{v}_N^{(n)})\}$
is bounded in $W^{1,2}(\Omega)$. Without loss of generality, we may assume that
$\{(\dot{v}_1^{(n)},\cdots,\dot{v}_N^{(n)})\}$ converges weakly in $W^{1,2}(\Omega)$
to an element $(\dot{v}_1,\cdots,\dot{v}_N)$. The compact embedding
\be\label{ex5}
W^{1,2}(\Omega)\to L^p(\Omega),\quad p\geq1,
\ee
implies that 
$
(\dot{v}_1^{(n)},\cdots,\dot{v}_N^{(n)})\rightarrow(\dot{v}_1,\cdots,\dot{v}_N)$
in $L^p(\Omega)$ $(p\geq1)$ as $n\rightarrow\infty$. In particular, $\underline{\dot{v}}_i=0$ $(i=1,\cdots,N)$.
In view of (\ref{3.14}) and (\ref{ex5}), we see that the functionals defined by the right-hand
sides of (\ref{3.31}) are continuous in $\dot{v}_i$ $(i=1,\cdots,N)$ with respect to
the weak topology of $W^{1,2}(\Omega)$. Therefore, $\underline{v}_i^{(n)}\rightarrow$ some $\underline{v}_i\in \mathbb{R}$
$(i=1,\cdots,N)$ as $n\rightarrow\infty$, where $\underline{v}_i$ is given in (\ref{3.31}).
In other words, ${\bf v}=(v_1,\cdots,v_N)=(\underline{v}_1+\dot{v}_1,\cdots,\underline{v}_N+\dot{v}_N)$
satisfies the constraints (\ref{3.26}) and solves the constrained minimization problem (\ref{3.27}).

\section{Solution on full plane}
\setcounter{equation}{0}

In this section, we prove the existence and uniqueness of the solution to
equations (\ref{3.1})--(\ref{3.2}) over $\mathbb{R}^2$ satisfying the natural boundary condition $u_i=2\ln v$ ($i=1,\cdots,N$) at infinity
as given in (\ref{bc}).
  Under the translation
$
u_i\mapsto u_i+2\ln v$ ($i=1,\cdots,N$) and the rescaling of the coefficient matrix, $v^2a_{ij}\mapsto a_{ij}$
($i,j=1,\cdots,N$),
  the equations (\ref{3.1}) become
   \bes
\Delta u_i=\sum_{j=1}^{N}{a}_{ij}(\e^{u_j}-1)+4\pi\sum_{s=1}^{n_i}\delta_{p_{i,s}}(x), \ \ \
i=1,\cdots,N,
 \lbl{4.3}\ees
where
   \bes
a_{ij}=\frac{v^2}{2N}\left({e^2}-\frac{g^2}{N}\right)+\delta_{ij}\frac{g^2 v^2}{2N},\ \
i,j=1,\cdots,N,
   \lbl{4.4}\ees
subject to the boundary condition 
    \bes
u_i\rightarrow 0,\ \ \ i=1,\cdots,N, \ \ \ {\rm as}\ |x|\rightarrow\infty.
   \lbl{4.5}\ees
   
As in \cite{JT,Ya1,Ya2}, we introduce the background function
    \bes
u_i^0(x)=-\sum_{s=1}^{n_i}\ln\left(1+\mu|x-p_{i,s}|^{-2}\right),\quad \mu>0,\quad i=1,\cdots,N.
 \lbl{4.7}\ees
Then we have
    \be
\begin{array}{lll}\Delta u_i^0&=&-h_i(x)+4\pi\sum_{s=1}^{n_i}\delta_{p_{i,s}}(x), \\
 h_i(x)&=&4\sum_{s=1}^{n_i}\frac{\mu}{(\mu+|x-p_{i,s}|^2)^2},\\
  \ \ \ i&=&1,\cdots,N.\end{array}
 \lbl{4.8}\ee
Let $u_i=u_i^0+U_i$, $i=1,\cdots,N$. Then the equations (\ref{4.3}) become
         \bes
\Delta U_i=\sum_{j=1}^{N}a_{ij}(\e^{u_j^0+U_j}-1)+h_i(x), \ \ \
i=1,\cdots,N.
         \lbl{4.9}\ees

Set
     \bess
&&\textbf{U}=(U_1,\cdots,U_N)^\tau,\ \ \ \textbf{G}=(\e^{u_1^0+U_1}-1,\cdots,\e^{u_N^0+U_N}-1)^\tau,\\
&& \textbf{H}=(h_1(x),\cdots,h_N(x))^\tau.
   \eess
Thus, the equations (\ref{4.9}) can be written in the vector form
    \bes
\Delta \textbf{U}=A\textbf{G}+\textbf{H},
    \lbl{4.10}\ees
where $A=(a_{ij})_{N\times N}$.

As before we use the transformation (\ref{3.7}) to change (\ref{4.10}) into
      \bes
\Delta v_i=t_{ii}\left(\e^{u_i^0+t_{ii}v_i+\sum_{k=1}^{i-1}\alpha_kv_k}-1\right)+\alpha_i\sum_{j=i+1}^N
\left(\e^{u_j^0+t_{jj}v_j+\sum_{k=1}^{j-1}\alpha_kv_k}-1\right)
+\sum_{j=1}^il_{ij}h_j, 
   \lbl{4.11}\ees
for $i=1,\cdots,N$.
It is direct to check that (\ref{4.11}) are the variational equations of the
energy functional
       \bes
I(\textbf{v})=\int_\Omega \sum_{i=1}^N\left\{\frac{1}{2}|\nabla v_i|^2+
\left(\e^{u_i^0+t_{ii}v_i+\sum_{k=1}^{i-1}\alpha_kv_k}-\e^{u_i^0}-\left[t_{ii} v_i+\sum_{k=1}^{i-1} \alpha_k v_k\right]\right)
+g_i v_i\right\}\,\ddd x,
\ees
where
\be 
g_i=\sum_{j=1}^il_{ij}h_j,\quad i=1,\cdots,N,
\ee
which may also be rewritten as
\bes
I({\bf v})&=&\sum_{i=1}^N\left\{\frac{1}{2}\|\nabla v_i\|_2^2+
\left(\e^{u_i^0},\e^{t_{ii}v_i+\sum_{k=1}^{i-1}\alpha_kv_k}-1-\left[t_{ii}v_i+\sum_{k=1}^{i-1}\alpha_kv_k\right]\right)_2\right.\nn\\
&&\quad +\left.
\left(\e^{u_i^0}-1,t_{ii}v_i+\sum_{k=1}^{i-1}\alpha_kv_k\right)_2+\left(g_i,v_i\right)_2\right\},\quad v_1,\cdots,v_N\in W^{1,2}(\mathbb{R}^2),
   \lbl{4.12}\ees
where $(\cdot,\cdot)_2$ and $\|\cdot\|_2$ denote the inner product and norm of $L^2(\mathbb{R}^2)$, respectively.

It is clear that the functional $I$ is a $C^1$-functional with respect to $\bf v$ and its Fr$\acute{e}$chet derivative satisfies
    \bes
\emph{D}I(\textbf{v})(\textbf{v})
=\sum_{i=1}^N\left\{\|\nabla v_i\|_2^2+
\left(t_{ii}v_i+\sum_{k=1}^{i-1}\alpha_kv_k,
\e^{u_i^0+t_{ii}v_i+\sum_{k=1}^{i-1}\alpha_kv_k}-1\right)_2
+(g_i,v_i)_2\right\}.
   \lbl{4.13}\ees

Set 
\be \label{x5.13}
w_i=t_{ii}v_i+\sum_{k=1}^{i-1}\alpha_kv_k,\quad i=1,\cdots,N,
\ee
${\bf w}=(w_1,\cdots,w_N)^\tau$, and ${\bf g}=(g_1,\cdots,g_N)^\tau$. Then ${\bf w}=T^\tau {\bf v}$ or ${\bf v}=L{\bf w}$. Thus
\be \label{x5.14}
\sum_{i=1}^N (g_i,v_i)_2=\int_{\mathbb{R}^2}{\bf g}^\tau{\bf v}\,\ddd x=\int_{\mathbb{R}^2}{\bf h}^\tau L^\tau L {\bf w}\,\ddd x=\sum_{i=1}^N(H_i,w_i)_2,
\ee
where ${\bf H}=(H_1,\cdots,H_N)^\tau=L^\tau L {\bf h}=A^{-1}{\bf h}$. Inserting (\ref{x5.13}) and (\ref{x5.14}) into (\ref{4.13}), we obtain
\bes
\emph{D}I(\textbf{v})(\textbf{v})
=\sum_{i=1}^N\left\{\|\nabla v_i\|_2^2+
\left(w_i,
\e^{u_i^0+w_i}-1+H_i\right)_2
\right\}.
   \lbl{x4.13}\ees
 
To
estimate the right-hand side of (\ref{x4.13}), we consider the quantity
   \bes
M(w)=\left(w,\e^{u_0+w}-1+H\right)_2,
   \lbl{4.14}\ees
where $w,u_0, H$ stand for one of the functions $w_i$, $u_i^0$, $H_i$, for $i=1,\cdots,N$, respectively.

As in \cite{JT}, we decompose $w$ into its positive and negative parts, $w=w_+-w_-$
with $w_+=\max\{w,0\}$ and $w_-=-\min\{w,0\}$. Then
$M(w)=M(w_+)+M(w_-)$. Using the inequality $\e^t-1\geq t$ ($t\in\mathbb{R}$),
we have
 $
\e^{u_0+w}-1\geq w+u_0,
$
which leads to
  \bes
M(w_+)\geq \|w_+\|^2_2+(w_+,u_0+H)_2\geq\frac12\|w_+\|^2_2-\frac12\|u_0+H\|^2_2,
    \lbl{4.16}\ees
where we have used the fact that $H,u_0\in L^2(\mathbb{R}^2)$.

On the other hand, using the inequality
$
1-\e^{-t}\geq {t}/{(1+t)}$ ($t\geq0$),
we can estimate $M(-w_-)$ from below as follows:
   \bes
M(-w_-)&=&\left(w_-,1-H-\e^{u_0}+\e^{u_0}[1-\e^{-w_-}]\right)_2
\nonumber\\
&\geq&\left(w_-,1-H-\e^{u_0}+\frac{w_-}{1+w_-}\e^{u_0}\right)_2\nonumber\\
&=&\left(\frac{w_-^2}{1+w_-},1-H\right)_2+\left(\frac{w_-}{1+w_-},1-H-\e^{u_0}\right)_2.
   \lbl{4.17}\ees
   
From the definition of the function $H$, we may choose $\mu>0$ large enough so that $H<1/2$ (say). Note also that $H,1-\e^{u_0}\in L^2(\mathbb{R}^2)$. Thus, we have
    \bes
\left|\left(\frac{w_-}{1+w_-},1-H-\e^{u_0}\right)_2\right|\leq \frac14\int_{\mathbb{R}^2}\frac{w_-^2}{1+w_-}\,\ddd x+\|1-H-\e^{u_0}\|^2_2.
   \lbl{4.18}\ees
Summarizing these facts, we see that (\ref{4.17}) enjoys the lower bound
    \bes
M(-w_-)\geq \frac{1}{4}\int_{\mathbb{R}^2}\frac{w_-^2}{1+w_-}\,\ddd x-C,
   \lbl{4.19}\ees
where and in the sequel, $C$ denotes a generic but irrelevant positive constant. 

From (\ref{4.16}) and (\ref{4.19}), we have
    \bes
M(w)\geq \frac{1}{4}\int_{\mathbb{R}^2}\frac{w^2}{1+|w|}\,\ddd x-C.
   \lbl{4.20}\ees

 Using (\ref{4.20}) in (\ref{x4.13}), we arrive at
    \bes\label{ex7}
\emph{D}I(\textbf{v})(\textbf{v})-\sum_{i=1}^N\|\nabla v_i\|_2^2
\geq\frac14\sum_{i=1}^N\int_{\mathbb{R}^2}\frac{w_i^2}{1+|w_i|}\,\ddd x-C.
\ees

Moreover, since the matrix $T$ is invertible and $\bf v$ and $\bf w$ are related through ${\bf w}=T^\tau{\bf v}$, we can find a positive constant $C_0$ such that 
\be \label{ex6}
\sum_{i=1}^N\|\nabla v_i\|_2^2\geq C_0\sum_{i=1}^N\|\nabla w_i\|^2_2.
\ee 

Inserting (\ref{ex6}) into (\ref{ex7}), we get
\be\label{ex8}
\emph{D}I(\textbf{v})(\textbf{v})\geq C_0\sum_{i=1}^N\|\nabla w_i\|_2^2
+\frac14\sum_{i=1}^N\int_{\mathbb{R}^2}\frac{w_i^2}{1+|w_i|}\,\ddd x-C.
\ee

We now recall the standard Gagliardo--Nirenberg--Sobolev inequality \cite{La,LU}
  \bes
\int_{\mathbb{R}^2}
f^4\,\ddd x\leq2\int_{\mathbb{R}^2}
f^2\,\ddd x\int_{\mathbb{R}^2}
|\nabla f|^2\,\ddd x, \ \ \ f\in W^{1,2}(\mathbb{R}^2).
   \lbl{4.24}\ees

Consequently, we have
\bea \label{3.21}
\left(\int_{\bfR^2}f^2\,\ddd x\right)^2&=&\left(\int_{\bfR^2}\frac{|f|}{1+|f|}(1+|f|)|f|\,\ddd x\right)^2\nn\\
&\leq&2\int_{\bfR^2}\frac{f^2}{(1+|f|)^2}\,\ddd x\int_{\bfR^2}(f^2+f^4)\,\ddd x\nn\\
&\leq&4\int_{\bfR^2}\frac{f^2}{(1+|f|)^2}\,\ddd x\int_{\bfR^2}f^2\,\ddd x\left(1+\int_{\bfR^2}|\nabla f|^2\,\ddd x\right)\nn\\
&\leq&\frac12\left(\int_{\bfR^2}f^2\,\ddd x\right)^2+C\left(1+\left[\int_{\bfR^2}\frac{f^2}{(1+|f|)^2}\,\ddd x\right]^4+\left[\int_{\bfR^2}|\nabla f|^2\,\ddd x\right]^4\right).\nn\\
&&
\eea
As a result of (\ref{3.21}), we have
\be \label{ex9}
\left(\int_{\bfR^2} f^2\,\ddd x\right)^{\frac12}\leq C\left(1+\int_{\bfR^2}|\nabla f|^2\,\ddd x+\int_{\bfR^2}\frac{f^2}{(1+|f|)^2}\,\ddd x\right).
\ee

 From (\ref{ex8}),
(\ref{ex9}), and the relation between $\bf v$ and $\bf w$, we may conclude with the coercive lower bound
       \bes
\emph{D}I(\textbf{v})(\textbf{v})\geq C_1
\left(\sum_{i=1}^N\|v_i\|_2+\sum_{i=1}^N\|\nabla v_i\|_2\right)
-C_2,\quad v_1,\cdots,v_N\in W^{1,2}(\bfR^2),
   \lbl{4.29}\ees
where $C_1,C_2$ are some constants.
In view of the estimate (\ref{4.29}), the existence of a critical point of the functional
$I$ in the space $W^{1,2}(\mathbb{R}^2)$ follows in a standard way. 

In fact, from (\ref{4.29}), we may
choose $R>0$ large enough such that
   \bes
\inf\left\{\emph{D}I(\textbf{v})(\textbf{v})\, |\,\textbf{v}=(v_1,\cdots,v_N)\in W^{1,2}(\mathbb{R}^2),
\ \|\textbf{v}\|_{W^{1,2}(\mathbb{R}^2)}=R\right\}\geq1.
    \lbl{4.30}\ees
    Consider the optimization problem
   \bes
\eta\equiv\inf\left\{I(\textbf{v})| \ \|\textbf{v}\|_{W^{1,2}(\mathbb{R}^2)}\leq R\right\}.
   \lbl{4.31}\ees
Let $\{\textbf{v}^{(n)}\}$ be a minimizing sequence of (\ref{4.31}). Without loss of
generality, we may assume that this sequence is also weakly convergent. Let
${\textbf{v}}$ be its weak limit. Thus, using the fact that the functional $I$ is weakly lower semi-continuous,
we have $I({\textbf{v}})\leq\eta$. Of course $\|{\textbf{v}}\|_{W^{1,2}(\mathbb{R}^2)}\leq R$
because norm is also weakly lower semi-continuous. Hence $I({\textbf{v}})=\eta$
and ${\textbf{v}}$ solves (\ref{4.31}). We show next that ${\textbf{v}}$ is a
critical point of the functional $I$. In fact, we only need to show that ${\textbf{v}}$ is
an interior point, or
 $
\|{\textbf{v}}\|_{W^{1,2}(\mathbb{R}^2)}< R.
  $
   For suppose otherwise that $\|{\textbf{v}}\|_{W^{1,2}(\mathbb{R}^2)}=R$. Then, in view of (\ref{4.30}),
 we have
     \bes
 \lim\limits_{t\rightarrow0}\frac{I({\textbf{v}}-t{\textbf{v}})-I({\textbf{v}})}{t}
 =\frac{\ddd}{\ddd t}I({\textbf{v}}-t{\textbf{v}})|_{t=0}=-\left(\emph{D}I({\textbf{v}})\right)({\textbf{v}})
\leq-1.
    \lbl{4.32}\ees
Therefore, when $t>0$ is sufficiently small, we see by virtue of (\ref{4.32}) that
$
I({\textbf{v}}-t{\textbf{v}})<I({\textbf{v}})=\eta.
$
However, since $\|{\textbf{v}}-t{\textbf{v}}\|_{W^{1,2}(\mathbb{R}^2)}=(1-t)R<R$, we arrive at
 a contradiction to the definition of ${\textbf{v}}$ or (\ref{4.31}). Thus $\bf v$ is a critical point of $I$.

 Finally, the strict convexity of $I$ says that $I$ can only have at most one critical point,
 so we have the conclusion that $I$ has exactly one critical point in $W^{1,2}(\mathbb{R}^2)$.
 Of course, this critical point is a solution of (\ref{4.3}), which must be
 smooth by virtue of the standard elliptic regularity theory.

The asymptotic behavior of the solution will be studied in the next section.

\section{Asymptotic behavior of planar solution}
\setcounter{equation}{0}

Let ${\bf v}=(v_1,\cdots,v_N)$ denote the solution of (\ref{4.3})  obtained in the previous section. Here we aim to establish the pointwise decay
properties for $\bf v$. 
Our tools are based on elliptic $L^p$-estimates and the maximum principle.

To proceed, we first recall the following embedding inequality  \cite{GT,JT,LU}
    \bes
\|f\|_{L^p(\mathbb{R}^2)}\leq \left(\pi\left[\frac{p}{2}-1\right]\right)^{\frac{p-2}{2p}}
\|f\|_{W^{1,2}(\mathbb{R}^2)},\ \ \ p>2.
   \lbl{5.1}\ees
From (\ref{5.1}), we may infer $\e^{f}-1\in L^2(\mathbb{R}^2)$ for $f\in W^{1,2}(\mathbb{R}^2)$.
To see this fact, we use a MacLaurin series expansion to get
    \bes
\left(\e^{f}-1\right)^2=f^2+\sum_{s=3}^\infty\frac{2^s-2}{s!}f^s.
   \lbl{5.2}\ees
By virtue of (\ref{5.1}) and (\ref{5.2}), we obtain
    \bes
\|\e^{f}-1\|^2_{L^2(\mathbb{R}^2)}\leq \|f\|_{L^2(\mathbb{R}^2)}^2+
\sum_{s=3}^\infty\frac{2^s-2}{s!}\left(\pi\frac{s-2}{2}\right)^{\frac{s-2}{2}}
\|f\|_{W^{1,2}(\mathbb{R}^2)}^s,
   \lbl{5.3}\ees
which confirms
our claim that $\e^f-1\in L^2(\bfR^2)$ because it is easily shown that the series on the right-hand side of (\ref{5.3}) is convergent. 

Now consider (\ref{4.11}). Since $v_i\in W^{1,2}(\mathbb{R}^2)$ ($i=1,\cdots,N$) and right-hand side of (\ref{4.11}) may be rewritten as
        \bes
&&t_{ii} \e^{u_i^0}\left(\e^{t_{ii}v_i+\sum_{k=1}^{i-1}\alpha_kv_k}-1\right)+
t_{ii}\left( \e^{u_i^0}-1\right)
+\alpha_i\sum_{j=i+1}^N\left\{  \e^{u_j^0}\left(\e^{t_{jj}v_j+\sum_{k=1}^{j-1}\alpha_kv_k}-1\right)+
\left( \e^{u_j^0}-1\right)\right\},\nn\\
&&i=1,\cdots,N,
   \lbl{5.4}\ees
we see that the right-hand side of each of the equations in (\ref{4.11})
belongs to $L^2(\mathbb{R}^2)$. Thus we may resort to
the standard elliptic $L^2$-estimates to deduce that $v_i\in W^{2,2}(\mathbb{R}^2)$ $(i=1,\cdots,N)$. In
particular, ${\bf v}(x)\rightarrow{\bf 0}$ as $|x|\rightarrow\infty$ because we are in two dimensions.

We can establish similar decay properties for $|\nabla v_i|$ $(i=1,\cdots,N)$. To this end, we rewrite the right-hand sides of (\ref{4.11})  as
        \bes
&&t_{ii}\left( \e^{u_i^0}-1\right)\e^{t_{ii}v_i+\sum_{k=1}^{i-1}\alpha_kv_k}+
t_{ii}\left( \e^{t_{ii}v_i+\sum_{k=1}^{i-1}\alpha_kv_k}-1\right)\nn\\
&&+\alpha_i\sum_{j=i+1}^N\left\{ \left( \e^{u_j^0}-1\right)\e^{t_{jj}v_j+\sum_{k=1}^{j-1}\alpha_kv_k}+
\left( \e^{t_{jj}v_j+\sum_{k=1}^{j-1}\alpha_kv_k}-1\right)\right\},\quad i=1,\cdots,N.
   \lbl{5.5}\ees
All these belong to the space $L^p(\mathbb{R}^2)$ for any $p>2$ due to the embedding
$W^{1,2}(\mathbb{R}^2)\to L^p(\mathbb{R}^2)$ and the definition of $u_i^0$.
Therefore 
elliptic $L^p$-estimates enable us to conclude that $v_i\in W^{2,p}(\mathbb{R}^2)$ $(i=1,\cdots,N; p>2)$. In particular,
we have $|\nabla v_i|\rightarrow0$ as $|x|\rightarrow\infty$, $i=1,\cdots,N$.

To obtain suitable exponential decay estimates for the solution,
it suffices to consider (\ref{4.3}) outside the disk $D_R=\{x\in\mathbb{R}^2\,| \,|x|<R\}$
where
         \bes
R>\max\{|p_{is}|\big| i=1,2,\cdots,N,\ s=1,2,\cdots,n_i\}.
   \lbl{5.6}\ees

For convenience, we write (\ref{4.3}) in $\mathbb{R}^2\setminus{D}_R$ in the form
   \bes
\Delta u_i=\sum_{j=1}^Na_{ij}u_j+\sum_{j=1}^Na_{ij}\left(\e^{u_j}-u_j-1\right), \ \
i=1,2,\cdots,N.
   \lbl{5.7}\ees
   
With the vectors ${\bf u}=(u_1,\cdots,u_N)^\tau$ and ${\bf V}=(\e^{u_1}-u_1-1,\cdots,\e^{u_N}-u_N-1)^\tau$, we may rewrite (\ref{5.7}) as
\be 
\Delta{\bf u}=A{\bf u}+A{\bf V}.
\ee 

Let $O$ be an $N\times N$ orthogonal matrix so that
   \bes
O^\tau AO={\rm diag}\{\lambda_1,\lambda_2,\cdots,\lambda_N\}\equiv\Lambda,
   \lbl{5.8}
            \ees
where $\lambda_1,\lambda_2,\cdots,\lambda_N$ are as given in (\ref{lambda}). Thus, in terms of the new variable vector
\bes
   \textbf{U}=(U_1,\cdots,U_N)^\tau=O^\tau{\bf u},
    \lbl{5.9}\ees
Substituting (\ref{5.9}) into (\ref{5.7}) and using (\ref{5.8}) and the behavior of
$\textbf{U}\rightarrow{\bf0}$ as $|x|\rightarrow\infty$, we arrive at
    \bes
\Delta \textbf{U}^2\geq 2\lambda_N\textbf{U}^2-a(x)\textbf{U}^2, \ \ \ x\in \mathbb{R}^2\setminus{D}_R,
   \lbl{5.11}\ees
where  $a(x)\rightarrow0$ as $|x|\rightarrow\infty$.
Consequently, for any $\varepsilon\in(0,1)$, we can find a suitably large $R_\varepsilon>R$
such that
      \bes
\Delta \textbf{U}^2\geq 2\lambda_N\left(1-\frac{\varepsilon}{2}\right)\textbf{U}^2, \ \ \ x\in \mathbb{R}^2\setminus{D}_{R_\varepsilon}.
   \lbl{5.12}\ees
Thus, using a suitable comparison function, the property ${\bf U}={\bf0}$ at infinity, and the maximum principle, we can obtain
a constant $C(\varepsilon)>0$ to achieve
   \bes
{\bf u}^2=\textbf{U}^2(x)\leq C(\varepsilon)\e^{-(1-\varepsilon)\sqrt{2\lambda_N}|x|},\quad |x|\geq R_\varepsilon.
   \lbl{5.13}\ees

We next derive some exponential decay estimates for $|\nabla u_i|$ ($i=1,\cdots,N$). For given $\ell=1,2$, we differentiate (\ref{5.7}) to obtain
    \bes
\Delta(\partial_\ell u_i)=\sum_{j=1}^Na_{ij}(x)\e^{u_j}(\partial_\ell u_j), \ \
i=1,\cdots,N.
   \lbl{5.14}\ees
Set
    ${\bf v}=(\partial_\ell u_1,\cdots,\partial_\ell u_N)^\tau$ and $E(x)={\rm diag}\{
\e^{u_1(x)},\cdots,\e^{u_N(x)}\}$.
Then  the system (\ref{5.14}) becomes
   \bes
\Delta{\textbf{v}}=A{\textbf{v}}+A(E(x)-I_N){\textbf{v}},
   \lbl{5.15}\ees
where $I_N$ is the $N\times N$ unit matrix. Consequently, we have
   \bes
\Delta{\textbf{v}}^2&\geq&2{\textbf{v}}^\tau\Delta{\textbf{v}}\nonumber\\
&=&2{\textbf{v}}^\tau A{\textbf{v}}+2{\textbf{v}}^\tau A(E(x)-I_N){\textbf{v}}\nonumber\\
&=&2\left(O^\tau{\textbf{v}}\right)^\tau\Lambda\left(O^\tau{\textbf{v}}\right)
+2{\textbf{v}}^\tau A(E(x)-I_N){\textbf{v}}\nonumber\\
&\geq&2\lambda_N{\textbf{v}}^2-a(x){\textbf{v}}^2, \ \ \ x\in \mathbb{R}^2\setminus {D}_R,
  \lbl{5.16}\ees
where $a(x)\to0$ as $|x|\to\infty$. Hence, as before, we obtain the estimate
\be 
\sum_{i=1}^N |\partial_\ell u_i|^2={\bf v}^2\leq C(\varepsilon) \e^{-(1-\varepsilon)\sqrt{2\lambda_N}|x|},\quad |x|\geq R_\varepsilon,\quad\ell=1,2.
\ee

\section{Consequences of asymptotic estimates}
\setcounter{equation}{0}

For the function $h_i$ given in (\ref{4.8}), we can directly compute to get
\be\label{7.1}
\int_{\mathbb{R}^2}h_i\,\ddd x=4\pi n_i,\quad i=1,\cdots,N.
\ee
Furthermore, for the solution $(U_1,\cdots,U_N)$ of the system (\ref{4.9}), we have $U_i=u_i-u^0_i$ where $u^0_i$ is defined by (\ref{4.7}) ($i=1,\cdots,N$) and
$(u_1,\cdots,u_N)$ is the unique solution of (\ref{4.3}) subject to the boundary condition (\ref{4.5}) whose derivatives  are seen to vanish at infinity exponentially fast. Thus
we infer that $|\nabla U_1|,\cdots,|\nabla U_N|$ all vanish at infinity at least as fast as $|x|^{-3}$. Consequently, we have
\be \label{7.2}
\int_{\mathbb{R}^2}\Delta U_i\,\ddd x=0,\quad i=1,\cdots,N.
\ee
Thus, integrating (\ref{4.9}) and applying (\ref{7.1}) and (\ref{7.2}) , we obtain the quantized integrals
\be\label{7.3}
\sum_{j=1}^N a_{ij} \int_{\mathbb{R}^2} (1-\e^{u_j})\,\ddd x= 4\pi n_i,\quad i=1,\cdots,N.
\ee

Now recall the relation between the functions $u_1,\cdots,u_N$ and the Higgs scalar fields $\phi_1,\cdots,\phi_N$ described in Section 2. Since
\be \label{7.4}
|\phi_i|^2 =v^2\e^{u_i},\quad i=1,\cdots,N,
\ee
we see that 
\be \label{7.5}
|\phi_i|^2-v^2=v^2(\e^{u_i}-1)=\mbox{O}(\e^{-(1-\varepsilon)\sqrt{2\lambda_N}|x|}),
\ee
if $|x|$ is large. Using (\ref{7.5}) in (\ref{2.10}), we conclude that the curvatures $B^i_{12}$ ($i=1,\cdots,N$) vanish at infinity exponentially fast at the same rate.
Besides, since the Higgs fields $\phi_1,\cdots,\phi_N$ and gauge fields $B^1_\ell,\cdots,B^N_\ell$  may be constructed from $u_1,\cdots,u_N$ following the expressions \cite{Ya2}
\bea 
\phi_i(z)&=&v \exp\left(\frac12 u_i(z)+\ii\sum_{s=1}^{n_i}\arg(z-p_{i,s})\right),\\
B^i_1(z)&=&-\mbox{Re}\{2\ii\overline{\partial}\ln \phi_i(z)\},\quad B^i_2(z)=-\mbox{Im}\{2\ii\overline{\partial}\ln \phi_i(z)\},
\eea
$i=1,\cdots,N$, we see that the covariant derivatives obey
\be 
\sum_{\ell=1}^2|(\partial_\ell -\ii B^i_\ell)\phi_i|^2=\frac{v^2}2\e^{u_i}|\nabla u_i|^2=\mbox{O}(\e^{-(1-\varepsilon)\sqrt{2\lambda_N}|x|}),\quad i=1,\cdots,N,
\ee
when $|x|$ is sufficiently large. Hence, inserting the value $\lambda_N=g^2/2N$, we obtain all the decay estimates stated in Theorem \ref{theorem2.2}.

Finally, applying (\ref{7.3}) and (\ref{7.4}) in (\ref{2.10}) and noting (\ref{4.4}), we obtain the quantized flux formulas stated in Theorem \ref{theorem2.3}
in the full plane case. 

In the situation of a doubly periodic domain, the same flux quantization conclusion follows simply from integrating (\ref{3.1}) and no further consideration is necessary.

\small{

}
 \end{document}